\documentclass{amsart}

\usepackage{latexsym}
\usepackage{amssymb}
\usepackage{graphicx}

\newtheorem{THM}{Theorem}
\newtheorem{Lemma}[THM]{Lemma}

\renewcommand{\Re}{\mathbb{R}}

\newcommand{\mH}{\mathbb{H}}

\newcommand{\Lap}{\Delta}

\renewcommand{\phi}{\varphi}
\newcommand{\mS}{\mathbb{S}}

\newcommand{\be}{\begin{equation}}
\newcommand{\ee}{\end{equation}}

\newcommand{\ra}{\rightarrow}
\newcommand{\ip}[2]{\left\langle#1,#2\right\rangle}
\newcommand{\lp}{\left(}
\newcommand{\rp}{\right)}
\newcommand{\lb}{\left[}
\newcommand{\rb}{\right]}

\DeclareMathOperator{\supp}{supp}

\DeclareMathOperator{\Ric}{Ric}
\DeclareMathOperator{\dive}{div}

\newcommand{\E}{\mathbb{E}}
\renewcommand{\P}{\mathbb{P}}
\newcommand{\QV}[1]{\left\langle#1\right\rangle}

\begin{document}

\title[Dirichlet problem at infinity]{Brownian motion and the Dirichlet problem at infinity on two-dimensional Cartan-Hadamard manifolds}
\author{Robert W.\ Neel}
\address{Department of Mathematics, Lehigh University, Bethlehem, PA, USA}
\begin{abstract}
After recalling the Dirichlet problem at infinity on a Cartan-Hadamard manifold, we describe what is known under various curvature assumptions and the difference between the two-dimensional and the higher-dimensional cases.  We discuss the probabilistic approach to the problem in terms of the asymptotic behavior of the angular component of Brownian motion.  Turning our attention to the two-dimensional case, we prove that the Dirichlet problem at infinity on a two-dimensional Cartan-Hadamard manifold is solvable under the curvature condition $K\leq (1+\epsilon)/(r^2 \log r)$, outside of a compact set, for some $\epsilon>0$ in polar coordinates around some pole.  This condition on the curvature is sharp, and improves upon the previously known case of quadratic curvature decay.  Finally, we briefly discuss the issues which arise in trying to extend this method to higher dimensions.
\end{abstract}
\thanks{The author gratefully acknowledges support from an NSF Postdoctoral Research Fellowship.}
\email{robert.neel@lehigh.edu}
\date{December 1, 2009}
\subjclass[2000]{58J32; Secondary 58J65 60H30}
\keywords{Dirichlet problem at infinity, Cartan-Hadamard manifold, Brownian motion}

\maketitle

\section{Introduction}

Given a Cartan-Hadamard manifold $M$, we can introduce the geometric boundary at infinity in the following way (see Section 3.1 of \cite{EltonBook} or Section 1 of \cite{EltonArticle} for an expanded treatment upon which the following is based).   Given a point $p$, the boundary at infinity $\mS_{\infty}(M)$ can be introduced by identifying it with the unit sphere in $T_pM$.  Further, if $(r,\theta)$ are polar coordinates around $p$, we can put a topology on $\widehat{M}=M\cup \mS_{\infty}(M)$ by saying that a sequence $x_n\in M$ converges to $\hat{\theta}\in \mS_{\infty}(M)$ if and only if $r(x_n)\ra\infty$ and $\theta(x_n)\ra \hat{\theta}$.  Both $\mS_{\infty}(M)$ and this topology can be shown to be independent of $p$, and the resulting topology on $\widehat{M}$ is called the cone topology.  It clearly makes $\widehat{M}$ compact.

For any continuous (real-valued) function $g$ on the boundary at infinity, the Dirichlet problem at infinity is the problem of finding a harmonic function $f$ on $M$ such that $f(x)\ra g(\hat{\theta})$ as $x\ra \hat{\theta}$ (in the cone topology).  For any Cartan-Hadamard manifold $M$, we say that the Dirichlet problem at infinity is solvable (on $M$) if it admits a unique solution for any $g$ as above.  A natural context for the Dirichlet problem at infinity is a Cartan-Hadamard manifold obeying a radial estimate on its sectional curvature, relative to some choice of pole.  We now discuss the relationship of the Dirichlet problem at infinity to other natural questions about the structure of harmonic functions on a Cartan-Hadamard manifold, and then indicate some known results.

First, we mention the conjecture of Greene and Wu that a Cartan-Hadamard manifold obeying the radial curvature estimate 
\begin{equation}\label{Eqn:GreeneWu}
K\leq -\frac{c}{r^2} \quad\text{when $r>R$}
\end{equation}
for some positive constants $c$ and $R$ admits a non-constant bounded harmonic function.  In this light, one effect of showing that the Dirichlet problem at infinity is solvable is to show that $M$ admits many non-constant bounded harmonic functions, and thus it provides a way of verifying this conjecture, at least on certain classes of manifolds.  Obviously, solvability of the Dirichlet problem at infinity does more than just produce non-constant bounded harmonic functions; it also relates them to the natural geometric notion of the boundary at infinity.  Continuing in this direction, one has the stronger question of identifying the Martin boundary (see Chapter 2 of \cite{SY} for background on the Martin boundary), and in particular, showing that the Martin boundary can be identified with the boundary at infinity for various classes of manifolds.  While solvability of the Dirichlet problem sheds some light on the relationship between the Martin boundary and the boundary at infinity, it stops short of showing they can be identified.

If the sectional curvature of $M$ is bounded above and below by negative constants, then Anderson and Schoen \cite{AS} have shown, using analytic methods, that the Martin boundary of $M$ is naturally homeomorphic to $\mS_{\infty}(M)$.  (Kifer \cite{Kifer} later developed a probabilistic version of the proof of Schoen and Anderson.)  As noted above, this is stronger than simply proving solvability of the Dirichlet problem, and essentially completely characterizes positive harmonic functions on such manifolds. 

Next, we observe that the situation is different depending on whether $M$ is two-dimensional, or of dimension greater than two.  If $M$ is two-dimensional, then the curvature bound
\begin{equation}\label{Eqn:SharpBound}
K \leq -\frac{1+\epsilon}{r^2 \log r} \quad\text{when $r>R$}
\end{equation}
for some positive constants $\epsilon$ and $R$ implies that $M$ is transient, and thus $M$ is conformally equivalent to the unit disk by uniformization.  This more than establishes the Greene-Wu conjecture in two dimensions.  Moreover, this bound is sharp, in the sense that a two-dimensional Cartan-Hadamard manifold satisfying the curvature bound
\[
K \geq -\frac{1}{r^2\log r} \quad\text{when $r>R$}
\]
for some positive $R$ has recurrent Brownian motion, or equivalently, is conformally equivalent to the plane.  In particular, such a manifold admits no non-constant positive harmonic functions.  Going further in the direction of understanding the relationship of harmonic functions to the boundary at infinity, Kendall and Hsu \cite{HsuKendall} prove, by studying the angular behavior of Brownian motion, that the Dirichlet problem at infinity is solvable on a two-dimensional Cartan-Hadamard manifold under the curvature bound $K\leq -c/r^2$ when $r>R$ for some positive constants $c$ and $R$ (in their paper, the authors note that their method cannot be pushed further to extend the result to the sharp curvature bound of Equation \eqref{Eqn:SharpBound}).

In dimensions three and higher, one also needs a lower curvature bound in order to show solvability of the Dirichlet problem at infinity.  In particular, Ancona \cite{Ancona} has produced a three-dimensional Cartan-Hadamard manifold with $K\leq -1$ everywhere (but not bounded from below) such that Brownian motion almost surely converges to a single point on the boundary at infinity, and a similar example where Brownian motion almost surely has every point on the boundary at infinity as an accumulation point.  Clearly, the Dirichlet problem at infinity is not solvable on either manifold.  However, it is important to note that both manifolds do admit non-constant bounded harmonic functions; these functions just don't relate to the boundary at infinity in the nice way that solutions to the Dirichlet problem do.  So Ancona's examples do not contradict the Greene-Wu conjecture, but they do indicate that comparing Brownian motion to geodesic rays (as probabilistic approaches to the Dirichlet problem at infinity do) can't resolve the conjecture in general.

In a positive direction, Hsu \cite{EltonArticle} has shown, also by studying the angular behavior of Brownian motion, that the Dirichlet problem at infinity is solvable in dimensions three and higher under the following upper and lower curvature bounds:
\[
K \leq -\frac{a(a-1)}{r^2}\quad\text{and}\quad \Ric \geq -r^{2b} \quad\text{when $r>R$}
\]
for some positive constants $R$, $a>2$, and $b<a-1$.  The same paper shows that a much more generous lower bound on the Ricci curvature can be allowed if $K$ is bounded from above by a constant, but here we're mostly interested in allowing the upper bound on $K$ to degenerate.

Finally, we mention the paper of Choi \cite{Choi}, which introduces a notion of convexity at infinity and uses it to study the Dirichlet problem at infinity by a type of Perron method.

As mentioned above, one natural approach to the Dirichlet problem at infinity is to study the angular behavior of Brownian motion on $M$.  In the following, we discuss the general structure of the problem, from a probabilistic point of view, and then revisit the case of two-dimensional Cartan-Hadamard manifolds, proving solvability of the Dirichlet problem at infinity with the sharp curvature estimate of Equation \eqref{Eqn:SharpBound}.  One nice feature of the argument is that it works directly with the semi-martingale decomposition of $\theta_t^2$, in contrast to previous probabilistic approaches.  Part of the appeal of this is that the martingale and bounded variation parts play different roles in the asymptotic behavior of Brownian motion (as explained in the next section), and this approach makes that explicit.  One drawback of the argument is that, at present, we have been unable to avoid an appeal to uniformization, although it is used only in a qualitative way in the proof of Lemma \ref{Lem:G2}.  (Then again, one can give a probabilistic proof, using elementary properties of the Green's function, that a transient surface of dimension two is conformally equivalent to the unit disk without too much difficulty, so maybe it's not that much of a drawback.)

\section{Probabilistic Background}

Again, we let $M$ be Cartan-Hadamard manifold, $p$ a point in $M$, and $(r,\theta)$ polar coordinates around $p$.  We will be interested in Brownian motion on $M$, denoted $B_t$, and we will use $r_t$ and $\theta_t$ to denote the composition of $r$ and $\theta$ with $B_t$.  Further, we use $\P^{\mu}$ and $\E^{\mu}$ to denote the probability and expectation with respect to Brownian motion on $M$ with initial condition given by a probability measure $\mu$.  In the case where $\mu=\delta_x$, we will abbreviate these as $\P^x$ and $\E^x$.  We let $\xi\in(0,\infty]$ be the explosion time of the Brownian motion.  Finally, we use $\QV{\theta}_t$ to denote the quadratic variation process of $\theta_t$.

In terms of Brownian motion, we have the following criterion for the solvability of the Dirichlet problem (see Proposition 6.1.1 of \cite{EltonBook})
\begin{THM}\label{THM:DPICriterion}
Let $M$ be as above.  Suppose that for every $x\in M$ we have that
\[
\P^x\lp \lim_{t\ra \xi} B_t \text{ exists (in the cone topology)}\rp =1 ;
\]
then we let $B_\xi$ denote this limit.  Further suppose that for any $\hat{\theta} \in \mS_{\infty}(M)$ and any neighborhood $N$ of $\hat{\theta}$ in $\mS_{\infty}(M)$, we have
\[
\lim_{x\ra \hat{\theta}}  \P^x\lp B_{\xi}\in N \rp = 1 .
\]
Then the Dirichlet problem at infinity is solvable on $M$.  Furthermore, for any continuous $g$ on $\mS_{\infty}(M)$, the (unique) solution to the Dirichlet problem with boundary function $g$ is given by $u(x)=\E^x \lb g(B_{\xi}) \rb$.
\end{THM}

Next, we discuss the special case when $M$ is radially symmetric (around $p$).  In this case, the metric can be written as
\[
ds^2 = dr^2 + J(r)^2 d\theta^2 
\]
where $d\theta^2$ is the usual metric on $\mS^{n-1}$ and $J(r)$, as the notation indicates, is a function only of $r$.  (Here ``J'' is meant to stand for ``Jacobi,'' since it is equal to the length of the obvious Jacobi fields.  Many authors use ``G'' instead, but we prefer to reserve this letter for the Green's function.)  Of course, $J$ is positive for positive $r$, $J(0)=0$, and $J^{\prime}=1$.  In this case, the radial component of Brownian motion satisfies the SDE (up to the explosion time $\xi$, which we generally won't say explicitly for future, similar equations)
\[
dr_t = dW_t +\frac{n-1}{2}\frac{\partial_r J}{J}\, dt ,
\]
where $W_t$ is a one-dimensional Brownian motion.  Further, the angular component of Brownian motion on $M$ is an independent (of $W_t$) time-changed Brownian motion on $\mS^{n-1}$, with time-change given by the integral of $1/J^2$ along Brownian paths.  Then convergence of the angular component of Brownian motion on $M$ essentially reduces to a one-dimensional problem.  Following these lines, March \cite{March} proved the following theorem (though stated somewhat differently).

\begin{THM}\label{THM:March}
Let $M$ be a Cartan-Hadamard manifold of dimension $n$, radially symmetric about a point $p$ (with polar coordinates $(r,\theta)$ around $p$).  Let $c_2=1$ and $c_n=1/2$ for $n\geq 3$.  Then if
\[
K \leq -\frac{c}{r^2\log r} \quad\text{when $r>R$} 
\]
for some $c>c_n$ and positive $R$, the Dirichlet problem at infinity is solvable (on $M$).  On the other hand, if
\[
K \geq -\frac{c}{r^2\log r} \quad\text{when $r>R$} 
\]
for some $c<c_n$ and positive $R$, there are no non-constant bounded harmonic functions on $M$.
\end{THM}

Now drop the assumption that $M$ is rotationally symmetric.  Then the angular component of Brownian motion will, in general, have a less symmetric martingale part and a drift term.  However, it is worthwhile to note that the same estimates used to prove the above theorem in the radially symmetric case still control the martingale part of $\theta_t$ in the general case, as a consequence of standard geometric comparison theorems.  Making this precise when $n\geq 3$ requires a bit of work and won't be needed in the following, so we won't do it.  However, we point out that, informally, it is the drift which is responsible for the bad behavior of the angular component of Brownian motion in the examples of Ancona mentioned above.

In the two-dimensional case, we have that $\theta$ is just the natural (real-valued) coordinate on $\mS^1$, and the metric can be written as
\[
ds^2 = dr^2 + J(r,\theta)^2 d\theta^2 ,
\]
where $J$ is the length of the obvious Jacobi field.  Then $B_t$ in polar coordinates satisfies the SDEs
\[
dr_t = dW_t +\frac{1}{2}\frac{\partial_r J}{J}\, dt \quad\text{and}\quad
d\theta_t = \frac{1}{J}\, d\widetilde{W}_t-\frac{1}{2}\frac{\partial_{\theta} J}{J^3}\, dt ,
\]
where $W_t$ and $\widetilde{W}_t$ are independent one-dimensional Brownian motions.  We now make the above claim about the martingale part of $\theta_t$ precise (in the two-dimensional case).

\begin{THM}\label{THM:2DQV}
Let $M$ be a Cartan-Hadamard surface, $p$ any point, and $(r,\theta)$ polar coordinates around $p$.  If the Gauss curvature $K(r,\theta)$ satisfies
\[
K \leq -\frac{1+\epsilon}{r^2 \log r} \quad\text{when $r>R$}
\]
for some positive constants $R$ and $\epsilon$, then for any $\delta>0$, $\P^{(r_0,\theta_0)}\lp\QV{\theta}_{\xi}>\delta\rp\ra 0$ as $r_0\ra\infty$, uniformly in $\theta_0$.
\end{THM}

\emph{Proof}:  Let $\tilde{K}(r)$ be smooth function on $[0,\infty)$, which is bounded from above by $0$ on $[0,R]$ and by $-(1+\epsilon)/(r^2 \log r)$ on $(R,\infty)$, and bounded from below by $\max_{\theta} K(r,\theta)$ (this is possible by the assumptions on $K$).  We think of $\tilde{K}$ as being the curvature of a radially symmetric comparison manifold.  We let $\tilde{J}(r)$ be the solution of the Jacobi equation determined by $\tilde{K}(r)$.  We know that $r_t$ satisfies the SDE (written in integral form, with initial value $r_0$)
\[
r_t = r_0+ W_t +\frac{1}{2}\int_0^t \frac{\partial_r J(r_s,\theta_s)}{J(r_s,\theta_s)}\, ds
\]
Then we let $\tilde{r}_t$ be the strong solution to
\[
\tilde{r}_t = r_0+W_t +\frac{1}{2}\int_0^t \frac{\partial_r \tilde{J}(\tilde{r}_s)}{\tilde{J}(\tilde{r}_s)}\, ds ,
\]
where we note that $\tilde{r}_t$ is driven by the same Brownian motion as $r_t$.  The Laplacian comparison theorem (see Theorem 3.4.2 of \cite{EltonBook}) implies that $\partial_r J/J \geq \partial_r \tilde{J}/\tilde{J}$, and thus a standard comparison result for SDEs (see Theorem 3.5.3 of \cite{EltonBook}) implies that $r_t \geq \tilde{r}_t$ for all $t\geq 0$ (up to the explosion time $\xi$), almost surely.  In addition, the Rauch comparison theorem implies that $J(r_1,\theta)\geq\tilde{J}(r_2)$ whenever $r_1\geq r_2$.  Since $\nabla \theta = 1/J$, it follows from these inequalities that the integral of $|\nabla \theta|^2$ along Brownian paths on $M$ is almost surely less than the integral of $1/\tilde{J}^2$ along $\tilde{r}_t$.

On the other hand, the same computations used to prove Theorem \ref{THM:March} show that the integral of $1/\tilde{J}^2$ along $\tilde{r}_t$ is almost surely finite, and that the probability that this integral exceeds any fixed positive level goes to zero as $r_0\ra\infty$.  (Indeed, the integral of $1/\tilde{J}^2$ along $\tilde{r}_t$ is easily realized as the quadratic variation of $\theta_t$ for a radially symmetric surface with curvature $\tilde{K}$.)  So we see that the integral of $|\nabla \theta|^2$ along Brownian paths on $M$ has these same properties, and the estimates hold independent of $\theta_0$, since the same comparison process $\tilde{r}_t$ can be used for any $\theta_0$.  This proves the theorem.  $\Box$

One thing to take away from the proof of above theorem is that an inequality for curvature implies inequalities for both the Laplacian of $r$ and the gradient of $\theta$, and these in turn give an inequality for the quadratic variation of the martingale part of $\theta_t$.  However, an inequality for curvature gives no control (at least pointwise) over the Laplacian of $\theta$, which poses a difficulty in our efforts to estimate the bounded variation part of $\theta_t$.  Indeed, in dimension three and higher, additional assumptions are needed for solvability of the Dirichlet problem.

\section{Dimension two}

In the two-dimensional case, however, it is possible to control the bounded variation part of $\theta_t$ without additional assumptions.  More specifically, choose any radial geodesic ray from $p$, and let $\theta_0$ be the angular coordinate corresponding to this ray.  We can define the truncated sectors, centered around $\theta_0$,
\[
S(\alpha,\beta) = \left\{ (r,\theta): r>\alpha, \theta\in (\theta_0-\beta,\theta_0+\beta) \right\}
\]
for $\alpha>0$ and $\beta\in(0,\pi)$.  We frequently write $S$, the dependence on $\alpha$ and $\beta$ (and $\theta_0$) being understood.  We will want to consider Brownian motion $B_t$ started at some probability measure $\mu$ supported in $S$, and stopped at the first exit time from $S(\alpha,\beta)$.  We let $G^{S}_{\mu}(r,\theta)$ be the associated Green's function on $S(\alpha,\beta)$.  When $\mu=\delta_x$, we write $G^S_x$.

To control the bounded variation part of $\theta_t$ in average, we essentially want to integrate $\Lap\theta/2$ against $G^S_{\mu}$.  However, just showing that this has small absolute value doesn't show that the bounded variation part isn't big (in absolute value), but rather that it's fairly ``symmetric.''  To get around this, we instead consider $(\theta_t-\theta_0)^2$.  The martingale part can be controlled just as for $\theta_t$ itself.  Further, if this semi-martingale has small expectation, then we can conclude that the bounded variation is small with high probability.  Thus, the focus of our argument is on estimating the expectation of $(\theta_t-\theta_0)^2$.  We begin with some preliminary lemmas.

\begin{Lemma}\label{Lem:T2}
For $M$ and $S(\alpha,\beta)$ as above, there exist positive constants $A$ and $C$ depending only on $R$ and $\epsilon$ such that
\[
\int_{S(\alpha,\beta)} |\nabla \theta|^2 \, d(\text{area}) \leq 2\beta C \frac{1}{\epsilon\lp\log\alpha\rp^{\epsilon}}
\]
for any $\alpha>A$.
\end{Lemma}
\emph{Proof:}  Recall that $|\nabla \theta|=1/J(r,\theta)$ and that, in polar coordinates, $d(\text{area}) = J(r,\theta) \,d\theta \,dr$.  Thus we have
\[
\int_{S(\alpha,\beta)} |\nabla \theta|^2 \, d(\text{area}) = 
\int_{\alpha}^{\infty} \int_{\theta_0-\beta}^{\theta_0+\beta} \frac{1}{J(r,\theta)} \, d\theta\, dr
\]
Now Rauch's comparison theorem and Equation \eqref{Eqn:SharpBound} imply (see \cite{Milnor} for the computation) that $J\geq r\lp \log r\rp^{1+\epsilon}/C$ for $r>A$, where $A$ and $C$ are positive constants depending only on $R$ and $\epsilon$.  Then, assuming that $\alpha>A$, we have
\[\begin{split}
\int_{S(\alpha,\beta)} |\nabla \theta|^2 \, d(\text{area}) & \leq
\int_{\alpha}^{\infty} \int_{\theta_0-\beta}^{\theta_0+\beta} \frac{C}{r\lp \log r\rp^{1+\epsilon}} \, d\theta\, dr \\
& = 2\beta C \lb \frac{-1}{\epsilon \lp\log r\rp^{\epsilon}}\rb_{\alpha}^{\infty} \\
&= 2\beta C \frac{1}{\epsilon\lp\log\alpha\rp^{\epsilon}} ,
\end{split}\]
and the lemma follows. $\Box$

\begin{Lemma}\label{Lem:G2}
Let $M$ be as above, and let $S=S(\alpha,\beta)$ be any truncated sector around $(r_0,\theta_0)$.  Then there exists $r_1>r_0$ and a probability measure $\mu$ such that the support of $\mu$ is $[r_0,r_1]\times\{\theta_0\}$ and $G^S_{\mu}$ is identically equal to one on $[r_0,r_1]\times\{\theta_0\}$.  Further, we have that
\[
\int_{S}\left|\nabla G^S_{\mu}\right|^2 \, d(\text{area}) = 2 ,
\]
where $\nabla G^S_{\mu}$ is interpreted in the weak sense on $\supp(\mu)$.
\end{Lemma}
\emph{Proof}:  We begin by showing the existence of $r_1$ and $\mu$.  For any $r_1$, let $p_{r_1}(x)$ be the function giving the probability that Brownian motion started at $x$ (and stopped upon exiting $S$) hits $[r_0,r_1]\times\{\theta_0\}$.  Then $p_{r_1}$ is harmonic on $S\setminus\lp [r_0,r_1]\times\{\theta_0\}\rp$, identically equal to one on $[r_0,r_1]\times\{\theta_0\}$, and goes to zero as $x$ approaches the boundary of $S$ or as $r(x)$ goes to infinity.  Thus $p_{r_1}$ is the Green's function corresponding to some measure $\nu_{r_1}$ with $\supp\lp \nu_{r_1}\rp = 
[r_0,r_1]\times\{\theta_0\}$.  Note that while $\nu_{r_1}$ has positive, finite mass, it may not have mass one.  Thus, in order to show the existence of $r_1$ and $\mu$ as indicated in the lemma, it is sufficient to show that for some choice of $r_1$, $\nu_{r_1}$ has mass one.

Note that the mass of $\nu_{r_1}$ increases as $r_1$ increases.  In particular, the mass of $\nu_{r_1}$ is just the capacity of $[r_0,r_1]\times\{\theta_0\}$ relative to $S$ (in the sense of potential theory), for which such monotonicity is well-known (see Section 4.3 of \cite{Grigoryan} for a summary of capacity on Riemannian manifolds).  It's relatively easy to see that the capacity of $[r_0,r_1]\times\{\theta_0\}$ is continuous as a function of $r_1$ and that it goes to zero as $r_1$ decreases to $r_0$.  Thus, we need to show that the capacity of $[r_0,r_1]\times\{\theta_0\}$ eventually reaches one as we increase $r_1$.  This follows from the stronger statement that the capacity increases without bound as $r_1$ increases.  To see this, note that capacity (in two dimensions) is invariant under conformal mapping, and that $S$ and $[r_0,\infty)\times\{\theta_0\}$ can be mapped conformally to $\mH^2$ and some smooth divergent curve on $\mH^2$.  (As mentioned in the introduction, this is the only point in our argument where conformal mapping techniques are used.)  Because $\mH^2$ is homogenous and isotropic, the Green's functions based at any two points coincide under some isometry, and also depend only on the distance form the base point.  In particular, we have uniform control over the decay of the Green's function based at any point.  Then it's not hard to see that on a divergent curve one can find a measure of any given finite mass, supported on the curve, such that the Green's function associated to this measure is everywhere less than or equal to one.  It follows that the capacity of $[r_0,r_1]\times\{\theta_0\}$ increases without bound as $r_1$ increases, and thus we have shown that $r_1$ and $\mu$ as described in the lemma exist.

If we now take $\mu$ to be as described in the lemma, it remains to compute the integral of $\left|\nabla G^S_{\mu}\right|^2$.  Note that $G^S_{\mu}$ is smooth on $S\setminus\supp(\mu)$.  Now choose $0<a<b<1$ such that the level sets $\{G^S_{\mu}=a\}$ and $\{G^S_{\mu}=b\}$ are smooth; note that this will be true for almost every $a$ and $b$ (with respect to Lebesgue measure on the reals) by Sard's theorem.  We let $n$ be the outward unit normal on the level sets of $G^S_{\mu}$, where ``outward'' means out of the corresponding super-level sets, so that $\nabla G^S_{\mu}$ is a negative multiple of $n$ on any smooth level set.  Then using the harmonicity of $G^S_{\mu}$, integration by parts gives
\[
\int_{\{ a\leq G^S_{\mu} \leq b \}} \left|\nabla G^S_{\mu}\right|^2 \, d(\text{area}) =
a \int_{\{G^S_{\mu}=a \}} \frac{\partial G^S_{\mu}}{\partial n} d(\text{length}) 
-b \int_{\{G^S_{\mu}=b \}} \frac{\partial G^S_{\mu}}{\partial n} d(\text{length})  .
\]
Next, note that the integral of $\partial G^S_{\mu}/\partial n$ over any smooth level set is $-2$ by the divergence theorem and the fact that $\mu$ has mass one (recall that $\Lap G^S_{\mu}/2 =\mu$ in the distributional sense).  Thus we have
\[
\int_{\{ a\leq G^S_{\mu} \leq b \}} \left|\nabla G^S_{\mu}\right|^2 \, d(\text{area}) =2(b-a) .
\]
Letting $a$ go to 0 and $b$ go to 1 (through values corresponding to smooth level sets) exhausts $S\setminus\supp(\mu)$ and shows that the integral of $\left|\nabla G^S_{\mu}\right|^2$ over $S\setminus\supp(\mu)$ is equal to two.

Finally, to complete the proof it suffices to show that $\nabla G^S_{\mu}$ can be extended to all of $S$ weakly by just letting it be zero on $\supp(\mu)$ (and thus $G^S_{\mu}$ is in the Sobolev space of functions on $S$ with a weak first derivative in $L^2$).  To do this, let $v$ be a smooth vector field, compactly supported on $S$.  Then integration by parts over $S\setminus\supp(\mu)$ gives (dropping the super- and sub-scripts for ease of notation)
\[\begin{split}
\int_{S\setminus\supp(\mu)} & G\dive(v) \, d(\text{area}) = \\
& \int_{\supp(\mu)} G\ip{v}{n} \, d(\text{length})
 - \int_{S\setminus\supp(\mu)} \ip{\nabla G}{v} \, d(\text{area}) ,
\end{split}\]
where we recall that $\supp(\mu)$ is just the curve $[r_0,r_1]\times\{\theta_0\}$, and where the integral over $\supp(\mu)$ is meant in the usual ``two-sided sense'' which comes from thinking of $\supp(\mu)$ as (part of) the boundary of $S\setminus\supp(\mu)$ with $n$ the outward unit normal.  Then the boundary integral over $\supp(\mu)$ is zero since the ``two sides'' cancel, which is just a consequence of $G$ and $v$ extending continuously across $\supp(\mu)$.  Also, because $\supp(\mu)$ has measure zero, we can extend the integral of $G\dive(v)$ to $S$ without changing its value, and we can also extend the integral of $\ip{\nabla G}{v}$ to all of $S$, after setting $\nabla G$ to be zero on $\supp(\mu)$, without changing its value.  We conclude that
\[
\int_{S} G\dive(v) \, d(\text{area}) = - \int_{S} \ip{\nabla G}{v} \, d(\text{area}) ,
\]
which establishes the desired weak differentiability of $G^S_{\mu}$.  $\Box$

We're now in a position to prove our main theorem.

\begin{THM}
Let $M$ be a Cartan-Hadamard manifold, and assume that, for some point $p$ with polar coordinates $(r,\theta)$ around $p$, the curvature satisfies
\[
K(r,\theta) \leq -\frac{1+\epsilon}{r^2 \log r} \quad\text{whenever $r>R$,}
\]
for some $\epsilon>0$ and $R>0$.  Then the Dirichlet problem at infinity is solvable on $M$.
\end{THM}

\emph{Proof:}  We choose any $\theta_0$ and a truncated sector $S(\alpha,\beta)$ around it.  We wish to understand the behavior of $(\theta-\theta_0)^2_t$, along Brownian motion stopped upon leaving $S(\alpha,\beta)$.  Without loss of generality, we assume for the moment that $\theta_0=0$ in order to lessen the notational burden.  Choose some $r_0>\alpha$ and let $\mu$ be a probability measure as described in Lemma \ref{Lem:G2}.  We see that $\theta^2_t$ satisfies the SDE
\[
d\lp \theta^2\rp_t = 2\theta \left|\nabla \theta\right|\, dW_t +\lp \theta\Lap\theta +\left|\nabla \theta\right|^2 \rp\, dt .
\]
Note that $|\theta|$ is bounded by $\beta$ on $S(\alpha,\beta)$.  Thus the quadratic variation of $\theta^2_t$, started from $\mu$, is bounded from above by $\beta$ times the quadratic variation of $\theta_t$, and so Theorem \ref{THM:2DQV} implies that the probability of $\QV{\theta^2}_{\xi}$ exceeding any fixed positive bound goes to zero as $r_0\ra\infty$.  (That our Brownian motion is started from $\mu$ rather than a point measure isn't a problem.  The quadratic variation starting from $\mu$ is just the $\mu$-average of the quadratic variation starting from every point in $\supp(\mu)$, all of which have $r\geq r_0$.)

To control the bounded variation part of $\theta^2_t$, we want to integrate $\theta\Lap\theta +\left|\nabla \theta\right|^2 $ against $G^S_{\mu}$.  The integral of $\left|\nabla \theta\right|^2 $ is just the quadratic variation of $\theta_t$, and so we know that the probability that it exceeds any fixed positive bound goes to zero as $r_0\ra\infty$, again by Theorem \ref{THM:2DQV}.  So we're left to estimate the integral of $\theta\Lap\theta $ against $G^S_{\mu}$.  We consider $S\cap\{r<\tilde{r} \}$ for some large $\tilde{r}$ in order to have a compact domain.  Then integration by parts gives
\[\begin{split}
\int_{S\cap\{r<\tilde{r} \}} & \lp\theta G^S_{\mu} \rp \Lap\theta \,d(\text{area}) =
\int_{\{r=\tilde{r}\}} \theta G_{\mu}^S \ip{\nabla\theta}{n} \,d(\text{area}) \\
& -\int_{S\cap\{r<\tilde{r} \}} G^S_{\mu} \left| \nabla\theta \right|^2  \,d(\text{area}) 
- \int_{S\cap\{r<\tilde{r} \}} \theta \ip{\nabla G^S_{\mu}}{\nabla\theta}  \,d(\text{area}) ,
\end{split}\]
where $n$ is the outward unit normal.  Now $\ip{\nabla\theta}{n}=0$, so the boundary integral is zero for all values of $\tilde{r}$.  If we let $\tilde{r}\ra\infty$, the integral of $G^S_{\mu} \left| \nabla\theta \right|^2$ just gives the quadratic variation of $\theta_t$ again.  So we know it's finite, and we have sufficient control of its size.  Further, as $\tilde{r}\ra\infty$, the absolute value of the $\theta \ip{\nabla G^S_{\mu}}{\nabla\theta}$ integral is bounded by $\beta\int_S |\ip{\nabla G^S_{\mu}}{\nabla\theta} |\, d(\text{area})$.  So we now need to estimate $\int_S |\ip{\nabla G^S_{\mu}}{\nabla\theta} |\, d(\text{area})$, since this will show that $\int_{S} \lp\theta G^S_{\mu} \rp \Lap\theta \, d(\text{area})$ is well-defined and give an estimate on its absolute value.

To do this, we simply use the Cauchy-Schwarz inequality.  By Lemmas \ref{Lem:T2} and \ref{Lem:G2}, the squared $L^2$-norms of $\nabla\theta$ and $\nabla G^S_{\mu}$ (over $S(\alpha,\beta)$) are bounded by
\[
2\beta C \frac{1}{\epsilon\lp\log\alpha\rp^{\epsilon}} \quad\text{and}\quad 2 ,
\]
respectively, for $\alpha>A$ sufficiently large (with $A$ and $C$ depending only on $R$ and $\epsilon$).  It follows that
\[
\int_{S(\alpha,\beta)}  |\ip{\nabla G^S_{\mu}}{\nabla\theta} |\, d(\text{area}) \leq
\int_{S(\alpha,\beta)}  |\nabla\theta||\nabla G^S_{\mu}|\, d(\text{area}) \leq \frac{2\sqrt{\beta C/\epsilon}}{\lp\log\alpha\rp^{\epsilon/2}}
\]
In particular, we can now conclude that the probability that the bounded variation part of $\theta^2_t$ ever exceeds any given positive level goes to zero as $\alpha\ra\infty$, which we can allow as $r_0\ra\infty$.

The previous statement applies to paths stopped upon leaving $S(\alpha,\beta)$, and our next goal is to extend it to unstopped paths.  In order for a path to hit the ``sides'' of $S$, that is, the $\theta=\pm\beta$ pieces of the boundary, $\theta_t^2$ would have to hit $\beta^2$.  We just showed that the probability of this happening can be made arbitrarily small by taking $\alpha$ to be large.  Now consider the $r=\alpha$ piece of the boundary of $S$.  Because $M$ is transient, it is possible to let both $\alpha$ and $r_0$ go to infinity in such a way that the probability of any path starting at $r\geq r_0$ hitting $r=\alpha$ goes to zero.  So we see that by taking $r_0$ and $\alpha$ large enough, the probability of Brownian motion started from the corresponding $\mu$ hitting the boundary of $S(\alpha,\beta)$ or having $\theta_t^2$ ever exceed any given positive level can be made arbitrarily small.  We conclude that, for Brownian motion started from $\mu$ and allowed to run all the way to $\xi$, the probability that $\theta^2_t$ ever exceeds any given positive level can be made arbitrarily small by choosing $r_0$ large enough.

Next, note that the probability of an event for Brownian motion stated from $\mu$ is just the $\mu$-average of the probability of the event for Brownian motion started from the points in $\supp(\mu)$.  Also, recall that our above reasoning applies to the ray corresponding to any $\theta_0$ (we just assumed $\theta_0=0$ for ease of notation).  So we see that, for any $\theta_0$, there is a sequence $r_{0,n}$ such that, for Brownian motion started from $(r_{0,n},\theta_0)$, the probability that $|\theta_t-\theta_0|$ ever exceeds any given positive level goes to zero as $n\ra\infty$ (obviously, an estimate on $(\theta_t-\theta_0)^2$ gives an estimate on $|\theta_t-\theta_0|$).

The idea now is to use these points to ``corral'' all other points on $M$.  We choose any $\phi\in[0,2\pi)$ and any small, positive constants $\delta$ and $\tilde{\delta}$.  Then we can find points $(\tilde{r}, \phi-2\delta)$ and $(\tilde{r}, \theta_0+2\delta)$ and some $\rho> \tilde{r}$ with the following properties: if we start Brownian motion from either $(\tilde{r}, \phi-2\delta)$ or $(\tilde{r}, \rho+2\delta)$ then the $\theta$-component never varies from its initial value by more than $\delta$ with probability at least $1-\tilde{\delta}$, and Brownian motion started at a point with its radial component greater than $\rho$ has probability at least $1-\tilde{\delta}$ of never having its radial component hit $\tilde{r}$.  We now consider any point $(r_0,\theta_0)$ with $r_0>\rho$ and $\theta_0\in (\phi-\delta,\phi+\delta)$ (in other words, we take a point in the truncated sector $S(\rho,\delta)$ around $\phi$).  We wish to show that the probability of $|\theta_t-\phi|$ ever exceeding $3\delta$ is small.  Consider any continuous path $\omega$ starting at $(r_0,\theta_0)$ and ending at a point with $\theta$-component $\phi+3\delta$ while its radial component always stays above $\tilde{r}$.  Then it follows (essentially just from the fact that curves divide surfaces, at least locally) that such an $\omega$ intersects at least the fraction $1-\tilde{\delta}$ of Brownian paths started from $(\tilde{r}, \rho+2\delta)$, by our assumptions on $(\tilde{r}, \phi+2\delta)$.

We use this in the following way.  Let $h(x)$ be the function on the sector $\tilde{S}=\{(r,\theta):r>0, \phi-3\delta<\theta<\phi+3\delta\}$ which gives the probability that a Brownian motion started at $x$  never leaves $\tilde{S}$.  Then $h$ is harmonic on $\tilde{S}$ and is at least $1-\delta$ at $(\tilde{r}, \rho+2\delta)$.  By stopping Brownian motion from $(\tilde{r}, \rho+2\delta)$ when it hits the image of $\omega$ and using the martingale property, we see that there is a point on $\omega$ where $h$ is at least $1-2\tilde{\delta}$.  The same is true for any path $\omega$ starting at $(r_0,\theta_0)$ and ending at a point with $\theta$-component $\phi-3\delta$ while its radial component always stays above $\tilde{r}$, by looking at Brownian paths from $(\tilde{r}, \rho+2\delta)$.  However, if we consider Brownian motion started from $(r_0,\theta_0)$, then once a path hits a point where $h\geq 1-2\tilde{\delta}$, it has probability at least $1-2\tilde{\delta}$ of never leaving $\tilde{S}$.  Combining this with the fact that $r_t$ stays above $\tilde{r}$ with probability at least $1-\tilde{\delta}$, a simple computation shows that the probability of a Brownian motion started at $(r_0,\theta_0)$, as described above, ever leaving the truncated sector $S(\tilde{r},3\delta)$ centered around $\phi$ is less than $3\tilde{\delta}$.  Note that this holds for any $(r_0,\theta_0)$ in the truncated sector $S(\rho,\delta)$ around $\phi$.

This result has the following consequence.  Choose any positive integer $n$, and decompose $\mS_{\infty}(M)$ into the union of $A_i =[2\pi i/n,2\pi (i+1)/n]$ for $i\in\{ 0,\ldots n-1 \}$.  Then for any small $\delta$, we can find $\rho$ (depending only on $M$, $n$ and $\delta$) such that, if $r_0\geq \rho$ and $\theta_0\in A_i$, then Brownian motion started at $(r_0,\theta_0)$ has $\theta_t\in A_{i-1}\cup A_i \cup A_{i+1}$ for all $t\in[0,\xi)$ with probability at least $1-\delta$ (here $i-1$ and $i+1$ should be understood mod $n$).  This follows just by applying the previous result to $n$ different values of $\phi$ (corresponding to the centers of the $A_i$) and choosing the other constants appropriately.

At this point, we're finally ready to verify the hypotheses of Theorem \ref{THM:DPICriterion}.  First, note that  $\P^x\lp \lim_{t\ra \xi} B_t \text{ exists}\rp$ is one for some $x$ if any only if it's one for every $x\in M$.  If this probability is not one everywhere, then there must be two disjoint subsets of $\mS_{\infty}(M)$, which we can take to be $A=[a,b]$ and $B=[c,d]$ for some $0\leq a<b<c<d<2\pi$, and a sequence of points $(r_n,\theta_n)$ such that $r_n\ra\infty$ and Brownian motion started from $(r_n,\theta_n)$ has probability at least $c>0$ of having accumulation points in both $A$ and $B$, for all $n$.  However, this clearly contradicts the observation in the previous paragraph.  Thus, $\P^x\lp \lim_{t\ra \xi} B_t \text{ exists}\rp=1$ for all $x$.  Once we know this, the other hypothesis of Theorem \ref{THM:DPICriterion} also follows immediately from the previous paragraph.  Thus we see that the Dirichlet problem at infinity is solvable on $M$.  $\Box$

We make a couple of technical comments about the proof.  We've already noted that considering $(\theta_t-\theta_0)^2$ instead of $\theta_t-\theta_0$ is minor point.  Considering $\mu$ and then ``corralling'' other points instead of just looking at $G_x^S$ directly is more of an annoyance.  The problem here is that, while $\{G_x^{\cdot}>1\}$ (relative to either $S$ or $M$) is compact, there doesn't seem to be any obvious way to control the angle it subtends and thus its contribution to the expected drift.  So instead we make $G_{\mu}^{\cdot}$ everywhere less than or equal to one, which gives the conclusion we want for $\mu$, and then we're left to show that it also holds for individual points.  These details aside, the heart of the argument is just integration by parts and the fact that the $L^2$-norm of $|\nabla\theta|$ can be estimated in dimension two.

\section{Future directions}

Given that the Dirichlet problem at infinity on a Cartan-Hadamard manifold of dimension two is solvable under a sharp curvature bound, the natural remaining problem in this direction is to show that the Martin boundary can be identified with the sphere at infinity (and the set of minimal positive harmonic functions).  One way to attempt this would be to refine the above approach, perhaps by incorporating an h-transform.

The other, probably more interesting, extension of the above results would be to adapt them to the higher-dimensional case, with the goal of improving the upper bound of \cite{EltonArticle} (in fact, one of the motivations for looking at the two-dimensional case is to prepare for the higher-dimensional case).  However, this is not straightforward, and we close by outlining the situation.  If we choose a ray $\gamma$, which corresponds to a point in $\mS_{\infty}(M)$, then we let $\theta(x)$ be the angle between $\gamma$ and the radial segment from $p$ to $x$.  This gives a coordinate near $\gamma$, and the main task is to show that $\theta_t$ has small bounded variation part, analogously to what was done above.  The idea would be to consider integrating $\Lap\theta$ against $G_x^S$.   For $G_x^S$, the situation is better than in the two-dimensional case, since the diameter of $\{G\geq 1\}$ is uniformly bounded by comparison with $\Re^3$, and thus there's no need to introduce a measure $\mu$ as above.  If we follow the same line of reasoning as above, integration by parts reduces the task to understanding $\int_{\{G\leq 1\}} \ip{\nabla\theta}{\nabla G}$.  The estimates used for the $L^2$-norm of $|\nabla G|$ over $\{G\leq 1\}$ still hold in higher dimensions.  However, in dimensions three and higher, $|\nabla \theta|$ no longer controls the volume element in polar coordinates.  Even for the model case of $\mH^3$, the $L^2$-norm of $|\nabla \theta|$ over a truncated sector is infinite for any reasonable way of defining $\theta$.  Thus there's no hope of a naive Cauchy-Schwarz estimate of the sort we used in two dimension working here.  Of course, one might hope for a more sophisticated way of controlling the integral of $\ip{\nabla\theta}{\nabla G_x^S}$.

\section*{Acknowledgment}

The author thanks Elton Hsu for introducing him to this problem.

\bibliographystyle{amsplain}

\def\cprime{$'$}
\providecommand{\bysame}{\leavevmode\hbox to3em{\hrulefill}\thinspace}
\providecommand{\MR}{\relax\ifhmode\unskip\space\fi MR }
\providecommand{\MRhref}[2]{\href{http://www.ams.org/mathscinet-getitem?mr=#1}{#2}}
\providecommand{\href}[2]{#2}

\end{document}